\documentclass{amsart}
\usepackage{latexsym}
\usepackage{amsmath, amssymb, amsfonts, amsxtra}
\usepackage{amsthm}

\author{Richard Oberlin}
\email{oberlin@math.wisc.edu}
\address{University of Wisconsin-Madison, Mathematics Department, 480 Lincoln Dr, Madison WI 53706}

\title{A recursive bound for a Kakeya-type maximal operator}
\date{}

\newcommand {\leb}{\mathcal{L}}
\newcommand {\gra}{\mathcal{G}}
\newcommand {\ort}{\mathcal{O}}
\newcommand {\rea}{\mathbb{R}}
\newcommand {\sph}{\mathbb{S}}
\newcommand {\ma}{\mathcal{M}}
\newcommand {\matilde}{\widetilde{\mathcal{M}}}
\newcommand {\proj}{\mathrm{proj}}
\newcommand {\spa}{\mathop{\mathrm{span}}}
\newcommand {\na}{\mathcal{N}}
\newcommand {\dis}{\mathrm{dist}}
\newcommand {\wtgra} {\widetilde{\gra}}
\newcommand {\Holder} {H\"{o}lder}
\newcommand {\kcrit} {k_{cr}}

\newcommand {\wtF}{\widetilde{F}}
\newcommand {\ovF} {\overline{F}}
\newcommand {\composed} {\circ}

\newtheorem{lem}{Lemma}[section]
\newtheorem{thm}{Theorem}[section]

\newtheorem{conj}{Conjecture}[section]
\newtheorem{prop}{Proposition}[section]

\begin{document}

\begin{abstract}
A $(d,k)$ set is a subset of $\rea^d$ containing a translate of
every $k$-dimensional plane. Bourgain showed that for 
$2^{k-1}+k \geq d$, every $(d,k)$ set has positive Lebesgue measure. 
We give an $L^p$ bound for the corresponding maximal operator.  
\end{abstract}

\maketitle

\section{Introduction}

A measurable set $E \subset \rea^d$ is said to be a $(d,k)$ set if it contains 
a translate of every $k$-dimensional plane in $\rea^d$. Once the definition
is given, the question of the minimum size of a $(d,k)$ set arises. This question 
has been extensively studied for the case $k=1$, the Kakeya sets. It is known 
that there exist Kakeya sets of measure zero, and these are called Besicovitch sets. It is 
conjectured that all Besicovitch sets have Hausdorff dimension $d$. For $k \geq 2$, it is 
conjectured that $(d,k)$ sets must have positive measure, i.e. that there are no $(d,k)$ 
Besicovitch sets. These size estimates are related to $L^p$ bounds on two
maximal operators which we define below. 

Let $G(d,k)$ denote the Grassmannian manifold 
of $k$-dimensional linear subspaces of $\rea^d$. For $L \in G(d,k)$ we define
\[
\na^k [f](L) = \sup_{x \in \rea^d} \int_{x + L} f(y) dy
\]
where we will only consider functions $f$ supported in the unit ball $B(0,1) 
\subset \rea^d$.

A limiting and rescaling argument shows that if $\na^k$ is bounded for some $p < \infty$ 
from $L^p(\rea^d)$ to 
$L^1(G(d,k))$, then $(d,k)$ sets must have positive measure. By testing $\na^k$ on the 
characteristic function of $B(0,\delta)$, $\chi_{B(0,\delta)}$, 
one sees that such a bound may only hold for $p \geq \frac{d}{k}$. 
For $ L $ in $ G(d,k) $ and $ a \in \rea^d $ define the $ \delta $ plate
centered at $a$, $ L_\delta(a) $, to be the $\delta$ neighborhood in $ \rea^d $
of the intersection of $ B(a,\frac{1}{2}) $ 
with $ L + a $. Fixing $L$, considering $\na^k \chi_{L_\delta(0)}$, and using the fact that 
the dimension of $G(d,k)$ is $k(d-k)$ we see that a bound into $L^q(G(d,k))$ can only hold for $q \leq kp$.
This leads to the following conjecture, where the case $k=1$ is excluded due to the existence of 
Besicovitch sets.
\begin{conj} \label{conjbdna}
For
$2 \leq k < d, p > \frac{d}{k}, 1 \leq q \leq k p$
\[
\|\na^k f \|_{L^q(G(d,k))} \lesssim  \|f\|_{L^p(\rea^d)}.
\] 
\end{conj}
It is also useful to consider a generalization of the Kakeya 
maximal operator, defined for $L \in G(d,k)$ by
\[
\ma^k_{\delta}[f](L) = \sup_{a \in \rea^d} \frac{1}{\leb^d(L_\delta(a))} 
\int_{L_\delta(a)} f(y) dy
\]
where $\leb^d$ denotes Lebesgue measure on $\rea^d$.
Using an argument analogous to that in Lemma 2.15 of \cite{bg}, one may see that 
a bound 
\begin{equation} \label{bound}
\|\ma^k_\delta f \|_{L^1(G(d,k))} \lesssim \delta^{\frac{-\alpha}{p}} \|f\|_{L^p(\rea^d)}
\end{equation}
where $ \alpha > 0 $ and $ p < \infty $, implies that the Hausdorff dimension of any $ (d,k) $ 
set is at least $ d - \alpha $.
Considering $\ma^k_\delta \chi_{B(0,\delta)}$ and $\ma^k_\delta \chi_{L_\delta(0)}$, we formulate
\begin{conj} \label{conjbdma}
For $k \geq 1, p < \frac{d}{k}, q \leq (d-k)p'$
\[
\|\ma^k_\delta f \|_{L^q(G(d,k))} \lesssim  \delta^{k-\frac{d}{p}} \|f\|_{L^p(\rea^d)}.
\] 
\end{conj}
In \cite{fa1} Falconer showed that $\na^k$ is bounded from $L^{\frac{d}{k}+ \epsilon}(\rea^d)$ 
to $L^1(G(d,k))$ for 
$k > \frac{d}{2}$. Later, in \cite{bg}, Bourgain used a Kakeya maximal operator bound
combined with an $L^2$ estimate of the $x$-ray transform to show that $\na^k$ is bounded 
from $L^p(\rea^d)$ to $L^p(G(d,k))$ for $(d,k,p) = (4,2,2 + \epsilon)$ and $(d,k,p) = (7,3,3 + \epsilon)$. 
He then showed, using a recursive metric entropy estimate, that   
for $ d \leq 2^{k-1} + k $, $\na^k$ is bounded for a large unspecified $p$. 
Substituting in the proof Katz and Tao's more recent bound for the Kakeya maximal operator in \cite{kt}
\begin{equation} \label{ktkakbound}
\| \ma^1_\delta f \|_{L^{n+\frac{3}{4}}(G(n,1))} \lesssim \delta^{-\frac{3(n-1)}{4n+3}}
 \|f\|_{L^{\frac{4n+3}{7}}(\rea^n)}
\end{equation}
one now sees that this holds for $k > \kcrit(d)$
where $\kcrit(d)$ solves $d = \frac{7}{6} 2^{\kcrit-1} + \kcrit$.

By \Holder 's inequality, the following holds for for any $k$-plate $L_\delta$ and positive $f$
\[
\int_{L_\delta}f \ dx \lesssim
\delta^{\frac{d-k}{r'}} 
\left( \int_{L^\perp} \left(\int_{L+y} f(x)\ d\leb^k(x) \right)^r d\leb^{d-k}(y)   \right)^{\frac{1}{r}}.
\]
Combining this with the $L^q(L^r)$ bounds on the $k$-plane transform proved by Christ in 
Theorem A of \cite{ch}, we see that Conjecture \ref{conjbdma} holds with 
$ p \leq \frac{d+1}{k+1} $. Except for a factor of $\delta^{-\epsilon}$, the same bound
for $ \ma^k_\delta $ was proven with $k=2$ by Alvarez in \cite{al} using a geometric-combinatorial
``bush''-type argument. More recently, also see \cite{mi2}. For dimension estimates of sets 
containing planes in directions corresponding to certain submanifolds of $G(4,2)$, see \cite{rogers}.

Our main result is the following.
\begin{thm} \label{boundonnak}
Suppose $4 \leq k < d$ and $\kcrit(d) < k$. Then 
\begin{equation} \label{boundonnakeq}
\|\na^k f\|_{L^{p}(G(d,k))} \lesssim 
\|f\|_{L^{p}(\rea^d)}
\end{equation}
for $f$ supported in the unit ball and $p > \frac{d-1}{2}.$
If, additionally, we have $k-j > \kcrit(d-j)$ for some integer $j$ in $[1,k-4]$, 
then we may take $p > \frac{d-3}{2+j}+1.$ 
\end{thm} 

The number $p = \frac{d-3}{2+j}+1$ is approximate, and may be slightly improved through 
careful numerology.  
We prove Theorem \ref{boundonnak} by combining a recursive bound of $\ma^k_\delta$ with
Bourgain's $L^2$ estimate. This recursive bound is based on Bourgain's 
metric entropy argument, but is carried out in a manner which is more efficient
for $L^p$ estimates. 
For $k \leq \kcrit(d)$ this method may be adapted to give the
following bound on $\ma^k_\delta$
\begin{thm} \label{maxoptwocorollary}
Suppose $2 \leq k \leq \kcrit(d)$. Then
\[
\| \ma^k_\delta f \|_{L^{\frac{d}{2}}(G(d,k))} \lesssim 
\delta^{-\frac{2}{d}(\frac{3(d-k)}{7(2^{k-2})} - 1 + \epsilon)} 
\|f\|_{L^{\frac{d}{2}}(\rea^d)}.
\]
\end{thm}
Finally, if $k+1 < \kcrit(d+1)$ then it is preferable not to use the $L^2$ bound, giving
\begin{thm} \label{maxoponecorollary}
For $2 \leq k$
\begin{equation} \label{maxoponecorollaryeq}
\| \ma^k_\delta f \|_{L^{d+1}(G(d,k))} \lesssim \delta^{-\frac{2}{d+1}(\frac{3(d-k)}{7(2^{k-1})}+ \epsilon)} 
\|f\|_{L^{\frac{d+1}{2}}(\rea^d)}.
\end{equation}
\end{thm}

From Theorems \ref{maxoptwocorollary} and \ref{maxoponecorollary}, we see that the 
Hausdorff dimension of any $(d,k)$ set is at least 
\[
\mathrm{min}(d ,d- \frac{3(d-k)}{7(2^{k-2})} + 1, d- \frac{3(d-k)}{7(2^{k-1})}).
\] 
It should be noted that the dimension estimate provided by only applying Theorem \ref{maxoponecorollary}
is also a direct consequence of the metric entropy estimate in \cite{bg}. 
However, to the best of the author's knowledge, 
it has not previously appeared in the literature, even without the improvement permitted by (\ref{ktkakbound}).

\subsection*{Acknowledgements}
I would like to thank my advisor Andreas Seeger for mathematical 
guidance and for his suggestion of the topics considered in this article.
I would also like to thank Daniel Oberlin for carefully 
reading several drafts.

\section{Preliminaries}

We start with the definition of the measure we will use on $G(d,k)$.
Fix any $L \in G(d,k)$. For a Borel subset $F$ of $G(d,k)$ let
\[
\gra^{(d,k)}(F)=\ort(\{\theta \in O(d): \theta(L) \in F\})
\]
where $\ort$ is normalized Haar measure of the orthogonal group on $\rea^d$, $O(d)$. Typically 
we will omit $d$ and $k$, denoting 
the measure by $\gra$. By the transitivity of the action of $O(d)$ 
on $G(d,k)$ and the invariance of $\ort$, it is clear that the definition is 
independent of the choice of $L$. Also note that $\gra$ is invariant under the 
action of $O(d)$. By the uniqueness of uniformly-distributed measures (see \cite{mt}, pages 44-53),
$\gra$ is the unique normalized Radon measure on $G(d,k)$ invariant under $O(d)$.

It will be necessary to use two alternate formulations of $\gra$. For each $ \xi $ in 
$ \sph^{d-1} $ let $ T_\xi:\xi^\perp \rightarrow \rea^{d-1} $ be an orthogonal linear 
transformation. Then $ T_\xi^{-1} $ identifies $ G(d-1,k-1) $ with the $ k-1 $ dimensional 
subspaces of $ \xi^\perp $. Now, define $T:\sph^{d-1}\times G(d-1,k-1) \rightarrow G(d,k)$ by
\[
T(\xi,L)=\spa(\xi,T_\xi^{-1}(L)).
\]
Choosing $T_\xi$ continuously on the upper and lower hemispheres of $ \sph^{d-1} $, 
$T^{-1}$ identifies the Borel subsets of $G(d,k)$ with the completion of the Borel
subsets of $\sph^{d-1} \times G(d-1,k-1)$. Under this identification, by uniqueness of 
rotation invariant measure, we have
\begin{equation} \label{graproduct}
\gra(F) = \sigma^{d-1} \times \gra^{(d-1,k-1)}(T^{-1}(F)).
\end{equation}
where $\sigma^{d-1}$ denotes normalized surface measure on the unit sphere.

It is also true that any invertible linear map $U: \rea^d \rightarrow \rea^d$ acts on
$G(d,k)$. We will need to know how $\gra$ varies under this action. Again using
the invariance of $\gra$, we observe that
\begin{equation} \label{lebgra}
\gra(F) = c \leb^{kd}
(\{(v_1, \ldots, v_k) : v_j \in B(0,1) \subset \rea^d \text{ and } \spa(v_1, \ldots, v_k) \in F\}).
\end{equation}
Using (\ref{lebgra}) and noting that, for $0 \neq r \in \rea$,
\[
\spa(v_1, \ldots, v_k) \in F \Leftrightarrow \spa(r v_1, \ldots, r v_k) \in F,
\]
we see that
\[
\frac{|\mathrm{det}(U)|^k}{\|U\|^{kd}} \gra(F) \leq \gra(U(F)) 
\leq |\mathrm{det}(U)|^k \|U^{-1}\|^{kd} \gra(F)
\]
where $\| \cdot \|$ denotes the operator norm of a linear map.
Since $|\mathrm{det}(U)| \leq \|U\|^d$ and $|\mathrm{det}(U)|=|\mathrm{det}(U^{-1})|^{-1}$, we have
\begin{equation} \label{graUF}
\left( \|U\| \cdot \|U^{-1}\| \right) ^{-kd} \gra(F) \leq \gra(U(F)) 
\leq \left( \|U\| \cdot \|U^{-1}\| \right)^{kd} \gra(F).
\end{equation}

\subsubsection*{Remark}

One should know that there have been two incorrect proofs published
on the subject of $(d,2)$ sets. The first, in \cite{fa2}, is well known and it
is of the claim that there are no Besicovitch $(d,2)$ sets for any $d$. The second, 
in \cite{mi1}, is of the claim that $(d,2)$ sets have Hausdorff dimension $d$ for
every $d$. Since it is quite recent, we will observe where the error is made. In 
the main construction, a $2$-plate $P^\delta$ is isolated which intersects a large number
of other $2$-plates $\{P^\delta_k \}$. Then a $\frac{\delta}{\rho}$ separated set
$\{e_i\} \subset \left( \sph^{d-1} \cap P^\perp \right)$ is chosen, and the set of 
$3$-plates $\{\Pi_i^{\tilde{C}\delta}\}$ is considered where 
each $\Pi_i^{\tilde{C}\delta}$ has the same center as $P^\delta$ and is in the 
direction $\spa(P,e_i)$. The aim is to show that each $P^\delta_k$ is contained in 
one of the $\Pi_i^{\tilde{C}\delta}$. However, it is only shown that for each 
$y \in P^\delta_k$ there is an $i$ so that $y \in \Pi_i^{\tilde{C}\delta}$ and hence 
\[
P^\delta_k  \subset \bigcup_i \Pi_i^{\tilde{C}\delta}.
\] 
The only assumption placed on the $P_k$ is that their distance from $ P $
is approximately $ \rho $ (where $ \rho \gg \delta $ ). 
For $ d \geq 4 $ if we let $P = \spa(v_1,v_2) $ and 
$P_k = \spa( \sqrt{1- \rho^2}v_1 + \rho v_3, \sqrt{1- \rho^2} v_2 + \rho v_4)$
where the $ v_j$ are orthonormal,
it can be seen that $P_k$  satisfies this assumption. However $P_k^\delta$ cannot be 
contained in any such $\Pi_i^{\tilde{C}\delta}$.

\section{A recursive maximal operator bound}

Our main argument is in the proof of Proposition \ref{maxopone} below. 

\begin{prop} \label{maxopone}
Suppose $k \geq 2$, $2 \leq p \leq d+1$, $p \leq r  \leq \frac{p(d-1)}{(p-2)}$, and $\frac{pr}{p+r} \leq q \leq r$. Then a bound for $\ma^{k-1}_\delta$ on $L^p(\rea^{d-1})$ of the form
\begin{equation} \label{kakbound}
\|\ma^{k-1}_\delta f\|_{L^q(G(d-1,k-1))} \lesssim \delta^{-\frac{\alpha}{p}} 
\|f\|_{{L^p}(\rea^{d-1})}
\end{equation}
implies the bound of $\ma^k_\delta$ on $L^{\tilde{p}}(\rea^d)$
\[
\|\ma^k_\delta f\|_{L^{\tilde{q}}(G(d,k))} \lesssim 
\delta^{-\frac{\tilde{\alpha}}{\tilde{p}}} 
\|f\|_{L^{\tilde{p}}(\rea^d)}
\]
with 
\[
\tilde{\alpha} = \frac{r \alpha}{p+r} + \epsilon,\ \ \tilde{p} = \frac{p(d-1) + 2r}{p+r},\  \text{\ and\ } \tilde{q} = \frac{\tilde{p} q(p+r)}{pr}.
\]
\end{prop} 

For our applications we will always take $q=p$.
It is then useful 
to note that the bound given by Proposition \ref{maxopone} is that which would result from interpolation 
between certain $L^2$ and 
$L^{d-1}$ bounds, namely:
\[
\frac{\tilde{\alpha}}{\tilde{p}} = \frac{\beta \alpha}{2} + \epsilon,\ \ \ \  
\frac{1}{\tilde{p}} = \frac{\beta}{2} + \frac{(1-\beta)}{d-1},\ \ \ \ 
\frac{1}{\tilde{q}} = \frac{\beta}{2} + \frac{(1-\beta)}{\infty} 
\]
where $\beta = \frac{2r}{p(d-1) + 2r}.$
For $\alpha < d-2k$
this is better than the bound given by interpolation between the case $r=p$ and the known 
(sharp) $L^2$ bound. However, $\tilde{p}$ still never seems to be optimal relative to $\tilde{\alpha}$ 
in the sense of Conjecture \ref{conjbdma}, which may be explained by the fact that we expect 
$\ma^k_\delta$ to be 
bounded independently of $\delta$ for $p>\frac{d}{k}$ rather than $p > d-1$.

The choice $r=p$ yields the greatest reduction of 
$\alpha$, giving $\tilde{\alpha}=\frac{1}{2} \alpha + \epsilon$.
However, this also gives a relatively large $\tilde{p}=\frac{d+1}{2}$. Alternately, 
choosing $r=\frac{p(d-1)}{(p-2)}$ gives a small reduction of $\alpha$ and a relatively large reduction
of $p$, with $\tilde{p} = \frac{p(d-1)}{d-2 + p-1}$. Observe that with this choice of $r$ and $m \leq d-2$ 
\begin{equation} \label{decreaseofp}
p = \frac{(d-1)-1}{m} + 1 \text{\ \ gives\ \ } \tilde{p} \leq \frac{d-1}{m+1} + 1.
\end{equation}  
Since interpolation with the $L^\infty$ bound does not affect $\tilde{\alpha}$, we may take 
$\tilde{p} = \frac{d-1}{m+1}+1$ for purposes of iteration.

To obtain Theorem \ref{maxoponecorollary}, we use $k-1$ applications of Proposition \ref{maxopone} with $r=p$. We start with the bound (\ref{ktkakbound}) for
$n=(d-k+1)$, except that we take $q_0 = p_0 = \frac{4n+3}{7}$ which is permitted by H\"{o}lder's inequality. This gives the value $\alpha_0 = \frac{3(n-1)}{7} + \epsilon = \frac{3(d-k)}{7} + \epsilon$. After one application of Proposition $\ref{maxopone}$, we have 
$p_1 = \frac{(n+1)+1}{2}$, $q_1 = (n+1)+1$, and $\alpha_1 = \frac{\alpha_0}{2} = \frac{3(d-k)}{2 \cdot 7} + \epsilon$. We use H\"{o}lder's inequality again, to take $q=p$, before another application of Proposition \ref{maxopone}. Continuing this process gives (\ref{maxoponecorollaryeq}).

To prove Proposition \ref{maxopone}, we will need the following lemma which gives a 
sort of parameterization of disjoint pairs of separated elements of a subset of an interval. 

\begin{lem} \label{subsets}
Suppose $\Gamma \subset [-1,1]$. Then for some integer $l$ with 
$\leb^1(\Gamma) \lesssim 2^{-l} \leq 1$ we have 
\[
\leb^1(S) \gtrsim \frac{1}{1 + |\log(\leb^1(\Gamma))|} \leb^{1}(\Gamma)
\]
where
\[
S=\{t \in [0,2^{-l}): |\Gamma \cap \{t + \mathbb{Z} 2^{-l}\}| \geq 2 \}. 
\]
\end{lem}

\begin{proof}

For each integer $l \geq -1$ let
\[
S_l = \{t \in [0,2^{-l}) : \sum_{j = 0}^{2^{l+1} - 1} \chi_{\Gamma}(t - 1 + j2^{-l}) \geq 2 \}
\]
and 
\[
M_l = \int_{S_l} \sum_{j = 0}^{2^{l+1} - 1} \chi_{\Gamma}(t - 1 + j2^{-l}) dt.
\]
Choose $l_0$ so that $\frac{\leb^1(\Gamma)}{4} < 2^{-l_0} \leq \frac{\leb^1(\Gamma)}{2}$. 
Note that 
\begin{equation} \label{subsetseq1} 
M_{l_0} \geq \frac{\leb^1(\Gamma)}{2}, 
\end{equation} and that 
\begin{equation} \label{subsetseq2} 
M_{-1} = 0.
\end{equation}
We want to observe that for some integer $l \in (-1,l_0],$ 
\begin{equation} \label{subsetseq3}
\leb^1(S_l) \geq \frac{\leb^1(\Gamma)}{16 l_0}.
\end{equation}
To see this, note that for every $l$, 
\begin{equation} \label{subsetseq4}
M_l - M_{l-1} \leq 2 \leb^1(S_l).
\end{equation}
Hence, if (\ref{subsetseq3}) does not hold for any $l \in (-1,l_0]$ then by
(\ref{subsetseq1}), (\ref{subsetseq4}), and induction we have
\[
M_l \geq \frac{\leb^1(\Gamma)}{2} - (l_0 - l)\frac{\leb^1(\Gamma)}{8 l_0}
\]
for every $l \in [-1,l_0]$. This is impossible by (\ref{subsetseq2}), proving (\ref{subsetseq3}). The lemma
follows since $l_0 \lesssim 1 + |\log(\leb^1(\Gamma))|$.
 
\end{proof}

In the next lemma we show that the maximal operator $\ma^k_\delta$ is local 
in the sense that we only need to prove bounds for functions supported in a ball. 

\begin{lem} \label{local}
For $q \geq p $ and $ r > 0 $ the bound 
$ \| \ma_\delta^k f \|_{L^q(G(d,k))} \leq C \|f\|_{L^p(\rea^d)} $
for all f supported in $ B(0,r) $ implies the bound
$ \| \ma_\delta^k f\|_{L^q(G(d,k))} \leq \widetilde{C} C \|f\|_{L^p(\rea^d)} $ 
for all $ f \in L^p(\rea^d) $ where $ \widetilde{C} $ is independent of $\delta$.
\end{lem}

\begin{proof}
Assume $ f $ is positive. Note that since the bound holds for functions supported 
in the ball of radius $ r $ centered at $ 0 $, it holds for functions supported in any ball
of radius $ r $. Pick a covering $ \{B(x_j,r)\}_{j=1}^{\infty} $ of $ \rea^d $ 
where each point in $ \rea^d $ is contained in only a finite number, say $ c_{d,r} $, of the 
balls. Then any plate $ L_\delta(a) $ touches at most 
$ \tilde{c}_{d,r} $ of the balls. So for any $ L $, 
\[
\ma_\delta^k f(L) \leq 
\tilde{c}_{d,r} \sup_j \ma_\delta^k \chi_{B(x_j,r)}f(L) \leq
\tilde{c}_{d,r} ( \sum_j ( \ma_\delta^k \chi_{B(x_j,r)}f(L) )^q )^{\frac{1}{q}}
\]
Interchanging $ L^q $ and $ l^q $ and using our bound:
\begin{eqnarray*}
\| \ma_\delta^k f(L) \|_{L^q(G(d,k))} \leq
\tilde{c}_{d,r} ( \sum_j \| \ma_\delta^k \chi_{B(x_j,r)}f\|_{L^q(G(d,k))}^q)^{\frac{1}{q}} \\  \leq
\tilde{c}_{d,r} C ( \sum_j \| f\|_{L^p(B(x_j,r))}^q)^{\frac{1}{q}}
\leq \tilde{c}_{d,r} C ( \sum_j \| f\|_{L^p(B(x_j,r))}^p)^{\frac{1}{p}}
\\ \leq \tilde{c}_{d,r} c_{d,r}^{\frac{1}{p}} C \| f \|_{L_p(\rea^d)}
\end{eqnarray*}
\end{proof}

\begin{proof}[Proof of Proposition \ref{maxopone}]
We will prove the restricted weak-type estimate for sets supported in $ B(0,1) $.
This will give the full estimate for functions supported in $ B(0,1) $ by 
interpolation. The general case then follows by Lemma \ref{local}, since $\tilde{q} \geq \tilde{p}$. 

We will only consider $ \delta \leq \frac{1}{2} $.
Let $ E \subset B(0,1) \subset \rea^d $. Fix $ 0 < \lambda \leq 1 $ and let  
\[
F = \{L \in G(d,k) : \ma_\delta^k [ \chi_E ](L) > \lambda \}.
\]
We need to show that
\begin{equation} \label{rwf}
\leb^d (E) \gtrsim \delta^{\frac{\alpha r}{p+r}+\epsilon} \lambda^{\tilde{p}} 
\gra(F)^{\frac{\tilde{p}}{\tilde{q}}}.
\end{equation}
By the trivial $L^1$ bound, there is a $c > 0$ so that (\ref{rwf})
is satisfied for $\lambda \lesssim \delta^c$. Thus, we may assume that 
$|\log(\lambda)| \lesssim |\log(\delta)|$. 

Instead of dealing directly with $F$, we will use its factorization via $T^{-1}$.
Let 
\[
\wtF = T^{-1}(F) \subset \left( \sph^{d-1} \times G(d-1,k-1) \right)
\] 
and $\wtgra = \sigma^{d-1} \times
\gra^{(d-1,k-1)}$. Then, by (\ref{graproduct}), we have $\wtgra(\wtF) = \gra(F)$.

Let $\{e_1, \ldots, e_d \}$ be an orthonormal basis of $\rea^d$.
For each integer $i \in [1,d]$ let $W_i = \{\xi \in \sph^{d-1}: |\langle \xi,e_i \rangle| \geq \frac{1}{\sqrt{d}} \}$.
Then $\sph^{d-1} = \bigcup_i W_i$, and thus for some $i$
\begin{equation} \label{gooddirection}
\int_{W_i} \int_{G(d-1,k-1)} \chi_{\wtF} dM d\xi \geq \frac{1}{d} \wtgra(\wtF) \gtrsim \gra(F).
\end{equation}
After renumbering assume that $i=d$.

Let $H=\spa(e_1,\ldots,e_{d-1})$ and define, for $\xi \in W_d$, the projection along $\xi$ onto $H$
\[
P_\xi(x) = x - \frac{\langle x,e_d \rangle}{\langle \xi,e_d \rangle} \xi.
\]
Henceforth, consider $G(d-1,k-1)$ as the set of $k-1$-planes in the particular copy $H$  of $\rea^{d-1}$.
We want to observe that if $L=\spa(\xi,M)$ where $\xi \in W_d$
and $M \in G(d-1,k-1)$, then for any $a \in \rea^d$ we have 
\begin{equation} \label{plateprojection}
P_\xi(L_\delta(a)) \subset cM_{\delta}(P_\xi(a))
\end{equation}
where $M_{\delta}(P_\xi(a)) \subset H$ is a $(k-1)$-plate  
and $c$ depends only on $d$.
To see this we first note that any point $l \in L_\delta(a)$ can be written 
\[
l = a + b\xi + m + w, 
\]
where $b \in \rea, m \in M, w \in L^\perp, |m| 
\lesssim \sqrt{d} $ and $|w| \leq \delta$. 
Then
\[
P_\xi l = P_\xi a + m + P_\xi w. 
\]
But since $\dis(\xi,H) \geq \frac{1}{\sqrt{d}}$ and $|w| \leq \delta$ it follows that
\[
|P_\xi w| = \left| w - \frac{\langle w,e_d \rangle}{\langle \xi,e_d \rangle}\xi \right| \leq 
\delta (1 + \sqrt{d}) \approx : c\delta.
\]
Thus $P_\xi(L_\delta(a))$ is contained in the $c\delta$ neighborhood, $cM_{\delta}(P_\xi(a))$,
of $P_\xi(a) + (M \cap B(0,c))$. 

For every $t \in \rea$ let $H_t=H + t e_d$ and $E_t = E \cap H_t$. Note that $P_\xi$ is an isometry from 
$H_t$ to $H$, giving
\begin{equation} \label{givingiso}
\leb^{d-1}\left(E_t \cap  L_\delta\left(a\right)\right) 
= \leb^{d-1}\left(P_\xi\left(E_t \cap L_\delta\left(a\right) \right) \cap cM_{\delta}
\left(P_\xi\left(a\right)\right) \right).
\end{equation}

The set $\wtF$ consists of pairs $(\xi,M)$ such that $\spa(\xi,T_\xi^{-1}(M)) \in F$.
However, considering (\ref{plateprojection}), we should be interested in pairs
$(\xi,M)$ such that $\spa(\xi,M) \in F$. 
We obtain a set of such pairs by letting 
\[
\ovF = \{(\xi, P_\xi \composed T^{-1}_\xi (\widetilde{M})): \xi \in W_d \text{\ and\ } (\xi,\widetilde{M}) \in \wtF \}.
\]
We will use our change of coordinates to estimate $\wtgra(\ovF)$.
Note that, by the orthogonality of $\xi$ and $\xi^{\perp}$, for $x \in H$
\[
|x| \leq |P_\xi \composed T_\xi^{-1}(x)| \leq (1 + \sqrt{d})|x|.
\] 
Then $\|(P_\xi \composed T_\xi^{-1})^{-1}\| \leq 1$ and $\|P_\xi \composed T_\xi^{-1}\| \leq (1 + \sqrt{d}).$
Thus, by (\ref{graUF}) and (\ref{gooddirection})
\begin{eqnarray} \label{grafbar}
\int_{W_d} \gra^{(d-1,k-1)} (\{M: (\xi,M) \in \ovF \}) d\xi
\\ \nonumber \gtrsim \int_{W_d} \gra^{(d-1,k-1)} (\{M: (\xi,M) \in \wtF \}) d\xi
\\  \nonumber \gtrsim \gra(F).
\end{eqnarray}

For $ \xi \in W_d$ and $s \neq t \in [-1,1]$, define the subset of $H$
\[
B_\xi^{s,t} = P_\xi(E_{s}) \cap P_\xi(E_{t}). 
\]
We will use the assumed maximal operator bound to estimate
\[
\left( \int_{W_d}  \leb^{d-1}(B_\xi^{s,t})^{\frac{r}{p}} d\xi \right)^{\frac{p}{r}}.
\]
This will provide us with an estimate of 
$\left( \leb^{d-1}(E_{s}) \leb^{d-1}(E_{t}) \right)^{\frac{p+r}{2r}}$, 
effectively reducing the exponent of $\delta$, as we will now explain.
Consider $E_s$ and $E_t$ as subsets of $H$ by orthogonal projection. Then 
\[
B_\xi^{s,t}= \left( E_{s} \cap \left(E_{t} + 
\frac{s-t}{\langle \xi, e_d \rangle} \proj_H (\xi) \right) \right) -
\frac{s}{\langle \xi, e_d \rangle} \proj_H (\xi)
\]
and so
\[
\leb^{d-1} (B_\xi^{s,t}) = \chi_{E_s} * \chi_{-E_t} \left( \frac{s-t}{\langle \xi, e_d \rangle} \proj_H (\xi) \right)
\] 
where we use $*$ to denote convolution in $\rea^{d-1}$.
Since $\dis(W_d,H) > \frac{1}{\sqrt{d}}$ 
\[
d\sigma^{d-1}(\xi) \lesssim |s-t|^{-(d-1)} d\leb^{d-1} \left(\frac{s-t}{\langle \xi, e_d \rangle} \proj_H(\xi) \right).
\] 
Thus by changing variables, Young's inequality, and the fact that $r \geq p$ 
\begin{eqnarray} \label{byyoungs}
\left( \int_{W_d} \left(\leb^{d-1}\left(B_\xi^{s,t}\right) \right)^{\frac{r}{p}} 
d\xi \right)^{\frac{p}{r}}
\lesssim
|s-t|^{\frac{-(d-1)p}{r}}
\left( \int_{\rea^{d-1}} (\chi_{E_s} * \chi_{-E_t} (x))^{\frac{r}{p}} dx \right)^{\frac{p}{r}} 
\\ \nonumber
\lesssim |s-t|^{\frac{-(d-1)p}{r}} (\leb^{d-1}(E_s) \leb^{d-1}(E_t))^{\frac{p+r}{2r}}.
\end{eqnarray}

We want to use our known maximal operator bound to estimate an average over $s$ and $t$ of the 
left hand side of (\ref{byyoungs}).
For each $x \in H, \xi \in W_d$ let 
\[
\Gamma_{\xi,x} = \{t: x \in P_\xi(E_t)\}.
\]
Then, if $(\xi,M) \in \ovF$ we have $L := \spa(\xi,M) \in F$ and hence for some $a_L \in \rea^d$,
\begin{eqnarray*}
\lambda \delta^{d-k} \lesssim 
\leb^{d}(L_\delta(a_L) \cap E) = 
\int_{-1}^{1} \int_{cM_\delta(P_\xi(a_L))} \chi_{P_\xi(E_t \cap L_\delta(a_L))} dx\ dt
\\
\leq \int_{cM_{\delta}(P_\xi(a_L))} \int_{-1}^1 \chi_{P_\xi(E_t)} dt\ dx
= \int_{cM_{\delta}(P_\xi(a_L))} \leb^1(\Gamma_{\xi,x}) dx
\end{eqnarray*}
where the first equality follows from (\ref{givingiso}). Thus, considering $\leb^1(\Gamma_{\xi,x})$ 
as a function of $x$, 
\[
\ma^{k-1}_{\delta} [\leb^1(\Gamma_{\xi,c \cdot})] (M) \gtrsim \lambda.
\]
Since $(\xi,M)$ was an arbitrary element of $\ovF$ and $r \geq q$ we now have by (\ref{grafbar})
\begin{equation} \label{maxofGamma}
\left( \int_{W_d} \left( \int_{G(d-1,k-1)} (\ma^{k-1}_\delta [\leb^1(\Gamma_{\xi,c \cdot})] (M))^q dM 
\right)^{\frac{r}{q}} d\xi \right)^{\frac{1}{r}}
\gtrsim \lambda \gra(F)^{\frac{1}{q}}.
\end{equation}

On the other hand, applying our assumed maximal operator bounds gives
\begin{eqnarray} \label{otoh}
\left( \int_{W_d} \left( \int_{G(d-1,k-1)} (\ma^{k-1}_\delta [\leb^1(\Gamma_{\xi,c \cdot})] (M))^q dM 
\right)^{\frac{r}{q}} d\xi \right)^{\frac{1}{r}}
\\ \nonumber \lesssim \delta^{-\frac{\alpha}{p}} 
\left( \int_{W_d} \left( \int_{\rea^{d-1}} \leb^1(\Gamma_{\xi,x})^p dx 
\right)^{\frac{r}{p}} d\xi \right)^{\frac{1}{r}}.
\end{eqnarray}
Let 
\[
Z = \{(\xi,x) \in W_d \times \rea^{d-1} : \leb^1(\Gamma_{\xi,x}) \gtrsim \frac{\lambda}{2} \}
\]
and note that (\ref{maxofGamma}) and (\ref{otoh}) still hold if we replace
$\leb^1(\Gamma_{\xi,x})$ by $\chi_Z \leb^1(\Gamma_{\xi,x})$. For each $(\xi,x) \in Z$,we may apply Lemma \ref{subsets} to $\Gamma_{\xi,x}$ obtaining an $l_{\xi,x}$ such that 
\begin{eqnarray*}
\lambda \lesssim 2^{-l_{\xi,x}} \leq 1 \\
\text{\ and\ } \leb^{1}(S^{l_{\xi,x}}_{\xi,x}) \gtrsim \frac{1}{1+|\log(\lambda)|}\leb^{1}(\Gamma_{\xi,x}) \\ 
\text{\ where \ } S^l_{\xi,x}=\{t \in [0,2^{-l}) : |\{t + \mathbb{Z} 2^{-l} \} \cap \Gamma_{\xi,x}| \geq 2\}. 
\end{eqnarray*}
Now, 
\[
\chi_Z \leb^1(\Gamma_{\xi,x}) = 
\chi_Z \leb^1(\Gamma_{\xi,x})
\sum_{i=1}^{C(1 + |\log{\lambda}|)} \chi_{\{i\}}(l_{\xi,x}) 
\]
and thus, combining (\ref{maxofGamma}) and (\ref{otoh}) we may choose $l_0$ so that 
\begin{eqnarray*}
\left( \int_{W_d} \left( \int_{\rea^{d-1}} 
 \chi_{\{l_0\}}(l_{\xi,x}) \chi_Z \leb^1(\Gamma_{\xi,x})^p dx \right)^{\frac{r}{p}} d\xi \right)^{\frac{1}{r}} 
\gtrsim \frac{1}{1 + |\log(\lambda)|} 
\delta^{\frac{\alpha}{p}} \lambda \gra(F)^{\frac{1}{q}}
\end{eqnarray*}
and hence
\begin{eqnarray} \label{beforemaxop}
\left(\int_{W_d} \left( \int_{\rea^{d-1}} \leb^1(S^{l_0}_{\xi,x}) ^p dx \right)^{\frac{r}{p}} d\xi
\right)^{\frac{1}{r}} \gtrsim 
 \frac{1}{(1 + |\log(\lambda)|)^2} \delta^{\frac{\alpha}{p}} \lambda \gra(F)^{\frac{1}{q}}.
\end{eqnarray}

Recalling the appropriate definitions, we see that
\[
S^{l_0}_{\xi,x} = \{t: x \in \bigcup_{i \neq j} B_\xi^{t+i2^{-l_0},t+j2^{-l_0}} \} 
\]
where $i$ and $j$ range over $\mathbb{Z} \cap [-2^{l_0},2^{l_0})$. This gives
\begin{equation} \label{untangled}
\leb^1(S^{l_0}_{\xi,x}) = \int_0^{2^{-l_0}} \sup_{i \neq j} \chi_{B_\xi^{t + i2^{-l_0},t+j2^{-l_0}}}(x)dt.
\end{equation}
Noting that the $L^r_{\sph^{d-1}}L^p_{\rea^{d}}L^1_{\rea}{}L^\infty_{\mathbb{Z}^2}$ norm is dominated 
by the $L^1_{\rea}{}L^p_{\mathbb{Z}^2}L^r_{\sph^{d-1}}L^p_{\rea^{d}}$ norm, we may combine
(\ref{beforemaxop}) and (\ref{untangled}), obtaining
\begin{eqnarray*}
\frac{1}{(1 + |\log(\lambda)|)^2} \delta^{\frac{\alpha}{p}} \lambda \gra(F)^{\frac{1}{q}}
\\ \lesssim
\int_0^{2^{-l_0}} \left( \sum_{i \neq j}  
\left( \int_{W_d} \leb^{d-1}(B_\xi^{t + i2^{-l_0},t+j2^{-l_0}})^{\frac{r}{p}} d\xi 
\right)^{\frac{p}{r}} \right)^{\frac{1}{p}} dt.
\end{eqnarray*}
Now, combining this with (\ref{byyoungs}), we have
\begin{eqnarray*}
2^{-l_0 \frac{(d-1)}{r}} \frac{1}{(1 + |\log(\lambda)|)^2} \delta^{\frac{\alpha}{p}} \lambda \gra(F)^{\frac{1}{q}}
\\ \lesssim
\int_0^{2^{-l_0}} \left( \sum_{i \neq j}  
\left( \leb^{d-1}(E_{t + i 2^{-l_0}}) \leb^{d-1}(E_{t + j 2^{-l_0}}) 
\right)^{\frac{p+r}{2r}}\right)^{\frac{1}{p}} dt.
\end{eqnarray*}
Finally, 
\begin{eqnarray*}
\int_0^{2^{-l_0}} \left( \sum_{i \neq j}  
\left( \leb^{d-1}(E_{t+i 2^{-l_0} }) \leb^{d-1}(E_{t + j 2^{-l_0}}) 
\right)^{\frac{p+r}{2r}}\right)^{\frac{1}{p}} dt
\\ \leq
\int_0^{2^{-l_0}} \left( \sum_{i}  
\leb^{d-1}(E_{t + i 2^{-l_0}})^{\frac{p+r}{2r}}\right)^{\frac{2}{p}} dt
\end{eqnarray*}
and by \Holder's inequality and the conditions $2 \leq p \leq r$
\[
\int_0^{2^{-l_0}} \left( \sum_{i}  
\leb^{d-1}(E_{t + i2^{-l_0}})^{\frac{p+r}{2r}}\right)^{\frac{2}{p}} dt
\leq 2^{-l_0(1 - \frac{2}{p})} \leb^{d}(E)^{\frac{p+r}{rp}}.
\]

Summarizing
\[
2^{-l_0 (\frac{(d-1)}{r}+\frac{2}{p}-1)} \frac{1}{(1 + |\log(\lambda)|)^2} \delta^{\frac{\alpha}{p}} \lambda \gra(F)^{\frac{1}{q}}
\lesssim
\leb^d(E)^{\frac{p+r}{rp}}.
\]
Since $2^{-l_0} \gtrsim \lambda$ and $\frac{d-1}{r} + \frac{2}{p} - 1 \geq 0$ we have 
\[
\lambda^{\frac{p(d-1) + 2r}{rp}} \delta^{\frac{\alpha}{p} + \epsilon} \gra(F)^{\frac{1}{q}}
\lesssim \leb^d(E)^{\frac{p+r}{rp}}
\]
or
\[
\delta^{\frac{\alpha r}{p+r} + \epsilon} \lambda^{\frac{p(d-1)+2r}{p+r}} \gra(F)^{\frac{pr}{q(p+r)}} \lesssim \leb^d(E).
\]
\end{proof}

\section{The $L^2$ method} \label{fourthsection}
 
Reducing $\alpha$ by a factor of two, as in Proposition \ref{maxopone}, is not a substantial gain
for small $\alpha$. The following proposition gives $\tilde{\alpha}=\alpha -1$ with 
$\alpha \geq 1$ and a bound for $\na^k$ with $\alpha < 1$. It is proved using Bourgain's technique 
from Propositions 3.3 and 3.20 of \cite{bg} in which he showed bounds for 
$\na^k$ with  $(d,k)=(4,2)$ and $(d,k)=(7,3)$. For completeness we will repeat the argument. 

\begin{prop} \label{bourgainrecursive}
Suppose $k,p \geq 2$ and that a bound for $\ma^{k-1}_\delta$ on $L^p(\rea^{d-1})$ of the form 
\begin{equation*} 
\|\ma^{k-1}_\delta f\|_{L^p(G(d-1,k-1))} \lesssim \delta^{-\frac{\alpha}{p}} 
\|f\|_{{L^p}(\rea^{d-1})}
\end{equation*}
is known. Then if $\alpha \geq 1$ we have the bound
\begin{equation} \label{alphageq1}
\|\ma^k_\delta f\|_{L^{p}(G(d,k))} \lesssim 
\delta^{-\frac{\alpha - 1 + \epsilon}{p}}\|f\|_{L^{p}(\rea^d)}
\end{equation}
for $f \in L^p(\rea^d)$.
If $\alpha < 1$ we have the bound
\begin{equation} \label{negativealpha}
\|\na^k f\|_{L^{p}(G(d,k))} \lesssim 
\|f\|_{L^{p}(\rea^d)}
\end{equation}
for $f \in L^{p}(\rea^d)$ supported in $B(0,1)$.
\end{prop}

To obtain Theorem $\ref{boundonnak}$, we start from an application of Theorem \ref{maxoponecorollary} with $k_0 = k-(2+j)$ and $d_0 = d-(2+j)$. After using H\"{o}lder's inequality on the left side, this gives
\begin{equation} \label{maxoponecorollaryresult}
\| \ma^{k_0}_\delta f \|_{L^{\frac{d_0+1}{2}}(G(d_0,k_0))} 
\lesssim \delta^{-\left(\frac{d_0+1}{2}\right)^{-1}\left(\frac{3(d-k)}{7\left(2^{k-(2+j)-1}\right)}+ \epsilon\right)} 
\|f\|_{L^{\frac{d_0+1}{2}}(\rea^{d_0})}.
\end{equation}
The condition $k-j > \kcrit(d-j)$ ensures that 
\[
\frac{3(d-k)}{7 2^{k - (2+j) - 1}} + \epsilon < 2,  
\]
and hence further reduction in $\alpha$ is unnecessary. Thus, with our $j$ ``spare'' iterations, we apply Proposition \ref{maxopone} with the maximum $r$ to give a reduction in $p$. Noting that $\frac{d_0+1}{2}$ satisfies the left equation in (\ref{decreaseofp}) with $m=2$, we start from (\ref{maxoponecorollaryresult}) to obtain after the first iteration
\[
\| \ma^{k_1}_\delta f \|_{L^{\frac{d_1-1}{3}+1}(G(d_1,k_1))} 
\lesssim \delta^{-\left(\frac{d_1 - 1}{3} + 1\right)^{-1}\left(\frac{3(d-k)}{7\left(2^{k-(2+j)-1}\right)}+ \epsilon\right)} 
\|f\|_{L^{\frac{d_1 - 1}{3} + 1}(\rea^{d_1})},
\]     
where $k_1 = k_0 + 1 = k-(2+(j-1))$, and 
$d_1 = d_0 + 1 = d-(2+(j-1))$. In fact, there is some additional improvement in $\alpha$ and $p$ which we ignore. After $j-1$ further iterations, we obtain
\[
\| \ma^{k_j}_\delta f \|_{L^{\frac{d_j-1}{j+2}+1}(G(d_j,k_j))} 
\lesssim \delta^{-\left(\frac{d_j - 1}{j+2} + 1\right)^{-1}\left(\frac{3(d-k)}{7\left(2^{k-(2+j)-1}\right)}+ \epsilon\right)} 
\|f\|_{L^{\frac{d_j - 1}{j+2} + 1}(\rea^{d_j})}, 
\]  
where $k_j = k-2$ and $d_j = d-2$. 
We then apply 
Proposition \ref{bourgainrecursive} twice, using (\ref{alphageq1}) the first time
and (\ref{negativealpha}) the second time, to obtain (\ref{boundonnakeq}).

Theorem $\ref{maxoptwocorollary}$ is obtained by instead applying Theorem \ref{maxoponecorollary} with $k_0 = k-1$ and $d_0 = d-1$, and then applying 
 (\ref{alphageq1}) from Proposition \ref{bourgainrecursive} once.

To prove Proposition \ref{bourgainrecursive} we will need an $L^2(L^2)$ estimate for the $x$-ray 
transform which utilizes cancellation. 
For every $ k > 0 $ let $ \phi^k $ be a positive Schwartz function on $ \rea^k $ such 
that $ \phi^k \geq 1 $ on $ B(0,\frac{3}{2}) $ and the Fourier transform, $ \hat{\phi^k} $, of $\phi^k$ has compact support. 
For $ \xi \in \sph^{d-1} $ and $ x \in \xi^\perp $ define
\[
\overline{f}_\xi(x) = \int \phi^1(t) f(x + t \xi)  dt.
\] 

\begin{lem} \label{radontransform}
Suppose $ \hat{f} \equiv 0 $ in $ B(0,R).$ Then
\[
\int_{\sph^{d-1}} \int_{\xi^\perp} |\overline{f}_\xi(x)|^2 dx d\xi \lesssim
\frac{1}{R} \|f\|^2_{L^2(\rea^d)}.
\]
\end{lem}

\begin{proof}
Choose $ N $ so that $ \hat{\phi^1} $ is supported in $ (-N,N) $.
Applying Plancherel's theorem to the partial Fourier transforms in the $ \xi $
and $ \xi^\perp $ directions, we have for every $\xi \in \sph^{d-1}$
\[
\int_{\xi^\perp} |\overline{f}_\xi(x)|^2 dx =
\int_{\xi^\perp} \left|\int_\rea \hat{\phi^1}(t) \hat{f}(\zeta + t \xi) dt\right|^2 d\zeta. 
\]
Considering the support of $ \hat{\phi^1} $ and using \Holder's inequality 
we have
\begin{eqnarray*}
\int_{\sph^{d-1}} \int_{\xi^\perp} \left|\int_\rea \hat{\phi^1}(t) \hat{f}(\zeta + t \xi) dt \right|^2 d\zeta\ d\xi 
\\ \nonumber \leq 2 N \|\hat{\phi^1}\|_{L^\infty}^2 
\int_{\sph^{d-1}} \int_{\rea^d} |\hat{f}(y)|^2\ |\chi_{[-N,N]}(\langle y,\xi \rangle )|^2 dy\ d\xi.
\end{eqnarray*}
Then for any $ y $
\begin{equation} \label{rotateplane}
\int_{\sph^{d-1}} \chi_{[-N,N]}(\langle y,\xi \rangle ) d\xi = 
\sigma^{d-1}\left(\left\{\xi : \dis(\xi, y^\perp) \leq \frac{N}{|y|} \right\}\right)
\lesssim \frac{N}{|y|}.
\end{equation}
Since $ |y| \geq R $ in the support of $ \hat{f} $, we are done.
\end{proof}

We will want to take advantage of the fact that the averaging operator $\ma^k_\delta$  
should tend to localize the Fourier transform. To this effect, we will define
a modified version of our maximal operator. For $ L \in G(d,k) $ let
\[
\pi_\delta^L(x) = \phi^k(\proj_L(x)) \delta^{-(d-k)} \phi^{(d-k)}\left(\proj_{L^\perp}\left( \frac{x}{\delta}\right) \right). 
\]
Now, define 
\[
\matilde_\delta^k [f](L) = \sup_{a \in \rea^d}\int_{\rea^d} \pi_\delta^L(a + x) f(x) dx.
\]
Immediately, we see that for all positive $ f, \ma_\delta^k [f] \lesssim \matilde_\delta^k [f] $.
We will see that the reverse inequality also holds. 

Let $ \varphi $ be a Schwartz function on $ \rea^d $ so that $ \widehat{\varphi} \equiv 1 $
on $ B(0,1) $ and $ \widehat{\varphi} $ is supported in $ B(0,2) $. For every $ R > 0 $ let
$ \varphi_R = R^d \varphi(R \cdot) $. 

\begin{lem} \label{switchdelta}
Suppose $ \hat{f} $ is supported in $ B(0,R) $. Then for any $k$-plane $ L \in G(d,k) $ and
$ a \in \rea^d $ we have
\begin{equation} \label{sdsecond}
\int_{a + (L \cap B(0,\frac{1}{2}))} |f(x)| dx \lesssim \ma_\frac{1}{R}^k[|f|](L)
\end{equation}
while for any $ \delta > 0 $ there is the estimate
\begin{equation} \label{sdthird}
\ma_{\delta}^k[|f|](L) \lesssim \ma_{\frac{1}{R}}^k[|f|](L).
\end{equation}
Also, without any assumptions on the support of $ \hat{f},$
\begin{equation} \label{sdfirst}
\matilde_\delta^k [|f|](L) \lesssim \ma_\delta^k [|f|](L).
\end{equation}
\end{lem}

\begin{proof}
The statement (\ref{sdthird}) follows from (\ref{sdsecond}) by averaging.
Inequality (\ref{sdfirst}) can be proved by the same method used in the proof
of (\ref{sdsecond}). So we will only prove 
(\ref{sdsecond}).

By our assumption on $ f $, $ f = f * \varphi_R $ so
\begin{eqnarray*}
\int_{a + (L\cap B(0,\frac{1}{2}))} |f(x)| dx =
\int_{a + (L\cap B(0,\frac{1}{2}))} |f*\varphi_R(x)| dx \\ \nonumber \leq
\int_{\rea^d} |\varphi_R(y)| \int_{a - y + (L\cap B(0,\frac{1}{2}))} |f(x)| dx\ dy.
\end{eqnarray*}
Let $ {e_1, \ldots ,e_d} $ be an orthonormal basis of $ \rea^d $ where 
$ L = \spa(\{e_1, \ldots ,e_k\}). $
For each  $ z \in \mathbb{Z}^d $ let 
$ b^R_z = (\frac{2}{\sqrt{d}R} z_1 e_1 , \ldots , \frac{2}{\sqrt{d}R} z_d e_d). $
Let $ Q^R_1 = L \cap B(0,\frac{1}{R}) $ and $ Q^R_2 = L^\perp \cap B(0,\frac{1}{R}). $

Then
\begin{eqnarray*}
\int_{\rea^d} |\varphi_R(y)|\int_{a - y + (L\cap B(0,\frac{1}{2}))} |f(x)| dx dy  \\
\leq \sum_{z \in \mathbb{Z}^d} \int_{b^R_z + Q^R_1 \times \ Q^R_2}
|\varphi_R(y)|\int_{a - y + (L\cap B(0,\frac{1}{2}))} |f(x)| dx dy \\
\lesssim \sum_{z \in \mathbb{Z}^d} \sup_{y \in b^R_z + Q^R_1 \times Q^R_2} |\varphi_R(y)| 
\int_{b^R_z+Q^R_1} \frac{1}{R^{d-k}} \ma^k_\frac{1}{R}[|f|](L)    dy' \\
\lesssim \ma^k_\frac{1}{R}[|f|](L) \sum_{z \in \mathbb{Z}^d} \frac{1}{R^d } 
\sup_{y \in b^R_z + Q^R_1 \times Q^R_2} |\varphi_R(y)|.  
\end{eqnarray*}

But 
\begin{equation} \label{schwartztailspl}
\sum_{z \in \mathbb{Z}^d} \frac{1}{R^d } \sup_{y \in b^R_z + Q^R_1 \times Q^R_2} |\varphi_R(y)| =
\sum_{z \in \mathbb{Z}^d} \sup_{y \in b^1_z + Q^1_1 \times Q^1_2} |\varphi_1(y)| 
\end{equation}
and the right-hand side of (\ref{schwartztailspl}) is controlled independently of $L$ and $R$ since
$\varphi_1$ is a Schwartz function. 
\end{proof}

\begin{proof}[Proof of Proposition \ref{bourgainrecursive}]
We will start by proving (\ref{alphageq1}). It suffices to consider the case when $ f $ is positive and 
bounded. By Lemma $\ref{local}$ we may also assume that $f$ is supported in 
$ B(0,1) $. Also we will only consider, say, $\delta \leq \frac{1}{2}$.

By (\ref{graproduct}), we need to show that
\begin{equation*}
\left(\int_{\sph^{d-1}} \int_{G(d-1,k-1)} \ma^k_\delta [f](\spa(\xi,T_\xi^{-1}M))^{p} dM d\xi 
\right)^{\frac{1}{p}}
\lesssim \delta^{- \left( \frac{\alpha - 1}{p} + \epsilon\right) } 
\|f\|_{L^{p}(\rea^d)}.
\end{equation*}
Note that, by our assumption on the support of $f$, 
\begin{equation} \label{notethat}
\ma^k_\delta [f](\spa(\xi,T_\xi^{-1}M)) \leq  
\ma_\delta^{k-1} [\overline{f}_\xi \composed T_\xi^{-1}](M).
\end{equation}

By a change of variables, the fact that $T^{-1}_\xi$ is orthogonal, and Plancherel's theorem in one dimension,
\begin{equation} \label{andplanchtheorem}
\widehat{\overline{f}_\xi \composed T_\xi^{-1}}(\zeta) =
\int_{\rea} \widehat{\phi^1}(t) \hat{f}(T_\xi^{-1} \zeta + t \xi) dt.
\end{equation}
Let $ g = f * \varphi_{\frac{c}{\delta}} $. 
Then by (\ref{andplanchtheorem}), the support of $\widehat{\phi^1}$, and our restriction on 
$\delta$, 
\begin{equation} \label{ftequal}
\widehat{\overline{g}_\xi \composed T_\xi^{-1}} \equiv 
\widehat{\overline{f}_\xi \composed T_\xi^{-1}} 
\text{\ on \ } B(0,\frac{\tilde{c}}{\delta}).
\end{equation} 
Hence, using (\ref{ftequal}) for the equality and Lemma \ref{switchdelta} for 
the last inequality
\begin{eqnarray} \label{usinglemmasd}
|\ma_\delta^{k-1} [\overline{f}_\xi \composed T_\xi^{-1}](M)| \lesssim 
|\matilde_\delta^{k-1} [\overline{f}_\xi \composed T_\xi^{-1}](M)| =
|\matilde_\delta^{k-1} [\overline{g}_\xi \composed T_\xi^{-1}](M)| 
\\ \nonumber \lesssim \ma_\delta^{k-1} [|\overline{g}_\xi| \composed T_\xi^{-1}](M). 
\end{eqnarray}

We will use the Littlewood-Paley decomposition of $g$.
Let $ \psi_0 = \varphi $
and for $ j > 0 $ let $ \psi_j = 2^{jd} \varphi(2^j \cdot) - 2^{(j-1)d} \varphi(2^{(j-1)} \cdot) $. 
Note that 
\[
\sum_{j=1}^\infty \widehat{\psi_j} \equiv 1
\]
and that for $ j > 0 $, $ \widehat{\psi_j} $ is supported in the annulus centered 
at $ 0 $ with radii $ 2^{k-1} $ and $ 2^{k+1} $.
For each $ j \geq 0 $, let $ g_j = g * \psi_j $. Then, considering the support of $\hat{g}$
\begin{equation} \label{sumofgjs}
g = \sum_{j=0}^{\log \frac{c}{\delta}} g_j.
\end{equation}
Now, by (\ref{usinglemmasd}) and (\ref{sumofgjs})
\begin{eqnarray} \label{eqnLpsecond}
\left(
\int_{\sph^{d-1}} \int_{G(d-1,k-1)} |\ma_\delta^{k-1} [\overline{f}_\xi \composed T_\xi^{-1}](M)|^{p} dM d\xi 
\right)^{\frac{1}{p}} 
\\ \nonumber \leq \sum_{j=0}^{\log \frac{c}{\delta}}
\left(
\int_{\sph^{d-1}} \int_{G(d-1,k-1)} \ma_\delta^{k-1} [|\overline{g_j}_\xi| \composed T_\xi^{-1}](M)^{p} dM d\xi
\right)^{\frac{1}{p}}.
\end{eqnarray}

Because each $\widehat{\overline{g_j}_\xi \composed T_\xi^{-1}}$ is supported in $B(0,2^{j+1})$,  
inequality (\ref{sdthird}) from Lemma \ref{switchdelta} allows us to apply our assumed bound  
with $\delta \approx 2^{-j}$ to give
\begin{eqnarray} \label{maxgj}
\int_{G(d-1,k-1)} \ma_\delta^{k-1} [|\overline{g_j}_\xi| \composed T_\xi^{-1}](M)^{p} dM  
\\ \nonumber
\lesssim \int_{G(d-1,k-1)} \ma_{2^{-j}}^{k-1} [|\overline{g_j}_\xi| \composed T_\xi^{-1}](M)^{p} dM
\\ \nonumber
\lesssim
2^{j \alpha }
\int_{\rea^{d-1}} |\overline{g_j}_\xi \composed T_\xi^{-1}|^{p} dx  \lesssim
2^{j \alpha} \|g_j\|_{L^{\infty}}^{p-2}\int_{\xi^\perp} |\overline{g_j}_\xi|^2 dx
\end{eqnarray} 
where, for the last inequality, we use the assumption that $p \geq 2$.
Because each $\hat{g_j}$ is identically zero on $B(0,2^{j-1})$,
integrating (\ref{maxgj}) and using Lemma \ref{radontransform} gives
\begin{equation*} \label{usinglemmagives}
\int_{\sph^{d-1}} \int_{G(d-1,k-1)} \ma_\delta^{k-1} [|\overline{g_j}_\xi| \composed T_\xi^{-1}](M)^{p} dM d\xi
\lesssim
2^{j\left(\alpha - 1 \right)} \|g_j\|_{L^{\infty}}^{p-2} \|g_j\|_{L^2}^2.
\end{equation*}
Thus
\begin{eqnarray} \label{rhslessthan}
\sum_{j=0}^{\log \frac{c}{\delta}}
\left(
\int_{\sph^{d-1}} \int_{G(d-1,k-1)} \ma_\delta^{k-1} [|\overline{g_j}_\xi| \composed T_\xi^{-1}](M)^{p} dM d\xi
\right)^{\frac{1}{p}}
\\ \nonumber\lesssim
\sum_{j=0}^{\log \frac{c}{\delta}} 
2^{\frac{j}{p}\left(\alpha -1 \right)} 
\|g_j\|_{L^{\infty}}^{1-\frac{2}{p}} 
\|g_j\|_{L^2}^{\frac{2}{p}}
\\ \nonumber\lesssim
\|f\|_{L^{\infty}}^{1-\frac{2}{p}} 
\|f\|_{L^2}^{\frac{2}{p}}
\sum_{j=0}^{\log \frac{c}{\delta}} 
\left(2^{\frac{\alpha - 1}{p}}\right)^j
\lesssim 
\|f\|_{L^{\infty}}^{1-\frac{2}{p}} 
\|f\|_{L^2}^{\frac{2}{p}}
\delta^{-\frac{\alpha - 1}{p}}
\end{eqnarray} 
Combining (\ref{notethat}) (\ref{eqnLpsecond}) and (\ref{rhslessthan}), we see that it only remains to show
\[
\|f\|_{L^\infty}^{1-\frac{2}{p}} \|f\|_{L^2}^{\frac{2}{p}} \lesssim \|f\|_{L_{p}}.
\]
This will hold under the additional assumption that $ f $ is a characteristic function. 
Sacrificing an $\epsilon$ in the exponent, this is sufficient by interpolation.

The proof of (\ref{negativealpha}) is identical except that we use (\ref{sdsecond}) instead of (\ref{sdthird}), and in (\ref{sumofgjs})
we must sum to $\infty$ instead of $\log(\frac{c}{d})$. This will converge in the end, by our 
assumption $\alpha < 1$. 
\end{proof}

\end{document}